\documentclass[12pt]{amsart}
\usepackage{amsmath}
\usepackage{amssymb,epsfig,amssymb,amsthm,epic,eepic,color,amsfonts}
\usepackage{hyperref}
\usepackage[all]{xy}
\usepackage{fontenc}  
\usepackage[utf8]{inputenc}

\newtheorem{thm}{Theorem}[section]

\newtheorem{prop}[thm]{Proposition}

\newtheorem{defn}[thm]{Definition}
\newtheorem{defns}[thm]{Definitions}

\newtheorem{cor}[thm]{Corollary}

\newtheorem{remarque}[thm]{Remark}
\newtheorem{remarques}[thm]{Remarks}

\newtheorem{rien}[thm]{}
\numberwithin{equation}{section}

\newcommand{\be}{\begin{enumerate}}
\newcommand{\ee}{\end{enumerate}}
\newcommand{\bi}{\begin{itemize}}
\newcommand{\ei}{\end{itemize}}

\def\R{\mathbb{R}}

\def\om{\omega}

\def\ga{\gamma}

\def\al{\alpha}
\def\be{\beta}
 
\def\de{\delta}

\def\vp{\varphi}

\def\la{\lambda}

\def\si{\sigma}
\def\Si{\Sigma}

\def\ep{\varepsilon}

\def\p{\partial}

\def\nd{\noindent}
\def\bull{${}$\hfill$\Box$\\}
\def\proof{\nd {\bf Proof.\ }}

\def\u{\underline}

\newcommand{\fl}[1]{{\color{red} #1}}

\begin{document}
\today
\vskip 1cm
\begin{center}
{\sc Conic singularities and immediate transversality
\vspace{1cm}

 Fran\c cois Laudenbach}
\end{center}

\title{}
\author{}
\address{Laboratoire de
math\'ematiques Jean Leray,  UMR 6629 du
CNRS, Facult\'e des Sciences et Techniques,
Universit\'e de Nantes, 2, rue de la
Houssini\`ere, F-44322 Nantes cedex 3,
France.}
\email{francois.laudenbach@univ-nantes.fr}



\begin{abstract}  I introduced the notion of immediate transversality in a recent paper on $A_\infty$-structures
on Morse complexes. In the present paper, I get that immediate transversality holds 
for every compact submanifold with $C^1$ conic singularities, not only in the setup of Morse theory. 
This is an occasion of speaking of this type of singularities.
\end{abstract}
\maketitle
\thispagestyle{empty}
\vskip 1cm

This paper is an expanded version of a talk given 
in the geometry seminar 
at the mathematical center IRMA of the University of Strasbourg (France) on November 8, 2021.

In the first section, we recall the notion of submanifold with $C^1$ conic singularities. 
That is the type of singularities which
are visible in the closure of invariant submanifolds of Morse gradients under some easy assumption of \emph{simplicity} about the gradient of a given Morse function.
We introduce the notion of \emph{tangent cone}, in general different from the considered cone, except when the latter
is \emph{linear}. In the considered setup, linearizability is shown to be a general property.

Section \ref{conic-sing} consists of 
 a detailed reminder of what I proved on this topic in the late nineties. I add a complete proof 
concerning the $CW$-structure one can derive from such a simple Morse gradient. 

  Section \ref{imm-trans} is devoted to the notion of \emph{immediate transversality by flow} that I introduced last year. 
Namely, given a compact submanifold $\Si$ with $C^1$ conic singularities in an ambient manifold $M$,
find a vector field $X$ on $M$ whose flow $X^t$ puts $\Si$ immediately transverse 
to $\Si$, that is, for every $t\in (0,\ep)$, one has $X^t(\Si)\pitchfork \Si$. I proved  this property holds in the
setup of Morse theory and is a key for $A_\infty$-structures on Morse complexes, unique up to homotopy.

In the present paper, we get a slight generalization:
such an immediate transversality holds for any compact submanifold
with $C^1$ singularities. It would be interesting to know what is the actual framework where immediate transversality 
holds.

I thank Athanase Papadopoulos for having given me the occasion of speaking of this topic.

\section{Generalities about  $C^1$ conic singularities} \label{gen-conic-sing}

In this section we state and prove some general facts relative to this type of singularities.

\begin{defn} \label{conic}
A subset $\Si$ in a smooth\footnote{Here and in the whole paper, smooth will stand for $C^\infty$ in contrast to $C^1$.} $n$-dimensional manifold $M$ is said to be \emph{a submanifold with $C^1$ conic singularities}
if $\Si$ is the disjoint union of smooth submanifolds $\Si_j$, possibly empty, with $j=\dim\Si_j$ (named the strata), 
which have the following mutual behaviour: there exists a family of tubular neighborhoods $N_j$ around $\Si_j$
such that:
\begin{enumerate}
\item Each stratum of $\Si$ is transverse to the sphere bundle $\p N_j$; that is denoted $\Si\pitchfork \p N_j$.
\item The intersection $\Si\cap N_j$ is a $C^1$ subbundle of the disc bundle $N_j$ over $\Si_j$. 
More precisely, there are an open 
covering $\{U_\al\}_{\al\in A}$ of $\Si_j$
and $C^1$ trivializations $\phi_\al$ of the pair $(N_j, \Si\cap N_j)$ over $U_\al$ which reads
$$\phi_\al: \left(U_\al\times (D^{n-j}, cS)\right)\mathop{\longrightarrow}^{\cong} \,(N_j, \Si\cap N_j)\vert_{U_\al}.
$$
 \end{enumerate}
Here, $S$ is itself a submanifold of the sphere $S^{n-j-1}=\p D^{n-j}$ with $C^1$ conic singularities and $cS$
is the radial cone based on $S$.  In the sphere bundle, the fiber 
$\Si\cap \p N_{j,x}$  over $x\in U_\al$ is carried back by $\phi_\al$ to a submanifold  $\{x\}\times S$. 

The union $\Si_0^k:= \sqcup_{j=0}^k\Si_j$ is called the $k$-{\rm skeleton} of $\Si$.
\end{defn} 

\begin{remarques}\label{rem-conic1}
${}$

{\rm  
\nd {\bf 1.} By item (2), 
the above definition makes sense by an induction on the ambient dimension.  
Moreover,
if $\Si_k$ meets $N_{j,x}$ at some point $y$ then, in a trivializing chart, the tangent space $T_y\Si_k$ contains
a local parallel copy of $T_x\Si_j$ and hence $k>j$.
 So, only strata of dimension higher than $j$ are visible in $N_j\smallsetminus \Si_j$.

\nd {\bf 2.} If the tube $N_j$ has a radius sufficiently small, the \emph{a priori} non-radial cone 
$\Si_x$ is transverse to the sphere $t\,\p N_{j,x}$ in 
$N_{j,x}$ for every $t\in (0,1]$. This follows from
 the Taylor formula at order 1  with integral remainder for $C^1$ maps, applied here to $\phi_\al: \{x\}\times D^{n-j}\to N_j$
 with $x\in U_\al$, that is:
  \begin{equation}\label{taylor}
  \phi_\al({tz}) -D\phi_\al(0) \cdot tz=\int_0^1\left[D\phi_\al(stz) -D\phi_\al(0)\right]\cdot tz\, ds.
 \end{equation}
 }
 \end{remarques}\bull
 \medskip

  What follows in the remainder of this 
  section is only related to what happens in a tube around one stratum. So, we use lighter notation:
 $B$ is a---for simplicity---compact $k$-dimensional smooth manifold, $E$ is a smooth $(n-k)$-disk bundle over $B$
 and $C\subset E$ is a $C^1$ conic subbundle of $E$ over $B$. This means that the zero-section of $E$ is a stratum of 
 $C$ and each fiber $C_x$ is $C^1$ diffeomorphic to a radial cone $cS$ over a smooth submanifold  $S$
 with $C^1$ conic singularities in the $(n-k-1)$-dimensional sphere.
 
 There are local trivializations $\phi_\al: U_\al\times(D^{n-k}, cS)\to (E,C)\vert_{U_\al} $ where $\{U_\al\}$ is a finite
 covering of $B$. The change of trivialization over $U_{\al\beta}:= U_\al\cap U_\beta$
 is $\phi_{\beta\al}:= \phi_\beta^{-1}\phi_\al$ restricted to the part of the bundles over $U_{\al\beta}$. 
 This collection of trivializations is known to fulfill some \emph{cocycle} equation:
 \begin{equation}\label{cocyclic}
 \phi_{\al\ga}\phi_{\ga\beta} \phi_{\beta\al}= Id\vert_{U_{\al\beta\ga}\times (D^{n-k},cS)}
 \end{equation}
 where $U_{\al\beta\ga}$ stands for the triple intersection $U_\al\cap U_\beta\cap U_\ga$. By the chain rule applied
 to the vertical derivative at $0$ in a fiber over $x\in U_{\al\beta\ga}$, one gets a linear cocycle
 \begin{equation}\label{lin-cocycle}
 D^v\phi_{\al\ga}(0)\cdot D^v\phi_{\ga\beta}(0)\cdot D^v\phi_{\beta\al}(0)= Id\vert_{U_{\al\beta\ga}\times\R^{n-k}}
 \end{equation}
 In the next definition, $T^v_0E$ stands for the vertical tangent bundle to $E$ along 
 the zero-section.
 
 \begin{defn} The tangent cone bundle to $C$, noted $\mathcal TC$, is the subbundle of $T^v_0E$ over $B$
 defined by the linear 1-cocycle  formed by the collection of $\{D^v\phi_{\al\beta}(0)\}$ where $\al$ and $\beta$ range over
 the set of indices of the covering.
 \end{defn}
 
 Since two systems of local trivialization of $(E,C)$ define cohomologous cocycles, the same is true
 for vertical derivatives at 0. Therefore, up to an isomorphism, $\mathcal TC$ does not depend on the chosen 
 trivializations. 
 \begin{defn}\label{lin}
One says that $C$ is a \emph{linear} conic subbundle if it has  linear local trivializations $\phi_\al$ 
or, equivalently, 
$C$  is equal to its tangent cone. One says that $C$ is \emph{linearizable} if it is isomorphic to a linear conic subbundle. 
\end{defn}
In particular, if $C$ is linear each fiber is a radial cone. 

\begin{center}\label{radial}
 \begin{figure}[h]
 \includegraphics[scale =.6] {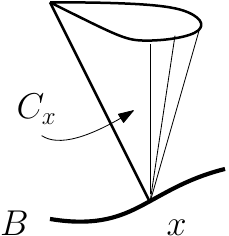}  
 \caption{\small  A schematic drawing for a one-parameter family of radial cones in 3-balls.}
\end{figure}
\end{center}

\begin{defn} One says that the pair $(E,C)$ is \emph{germinally linearizable} 
if there exists a disc subbundle $E'\subset E$
such that the sphere bundle $\p E'$ is transverse to $C$ and $C\cap E'$ is linearizable.\\
\end{defn}

\begin{prop} \label{linearization}
Every $C^1$ conic subbundle $C$ of a disc bundle $E$ is germinally linearizable.
\end{prop}

\proof By Remark \ref{rem-conic1}.{\bf 2}, there is a disc subbundle $E'\subset E$ such that for every 
$t\in(0,1]$ the sphere bundle $t\,\p E'$ is transverse to $C$. Set $C':= C\cap E'$.
Now, we apply to $C'$ a variation of the famous Alexander trick,  namely, 
for every $t\in (0,1]$
\begin{equation}
C'_t:= \frac 1t \cdot\left(C'\cap t\, E'\right)
\end{equation} 
\begin{center}\label{alexander}
 \begin{figure}[h]
 \includegraphics[scale =.6] {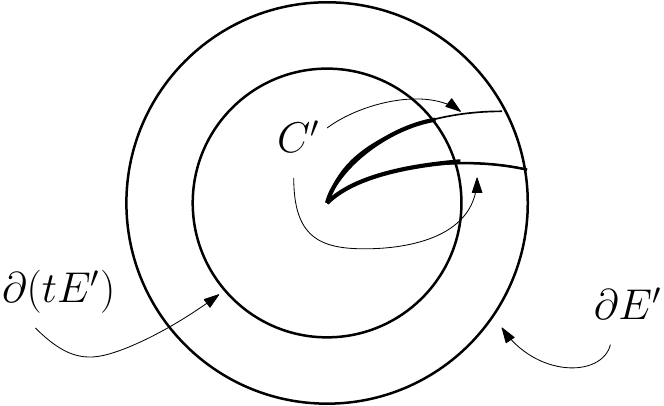}  
 \caption{\small  Cone $C'\cap (t\,E')$ is drawn with thick lines.}
\end{figure}
\end{center}
makes sense and defines some conic subbundle of $E'$. Its limit when $t$ goes to zero is the tangent cone bundle
$\mathcal TC\cap E'$ which follows from formula (\ref{taylor}).
This fiberwise isotopy of $C^1$ conic subbundles is induced by an ambient fiberwise isotopy of $E'$. Finally,
from this isotopy one derives an isomorphism to a linear conic bundle. For writing smooth formulas, a solution would be 
to choose a fiberwise contracting  vector field which is tangent to each stratum of $C'$ and uniquely integrable. 
This is allowed by the nature of  singularities. \bull

 Here are a few more generalities which allow us to build 
 new examples of submanifolds with $C^1$ conic singularities.
 
 \begin{defn} \label{strat-transv}
 Given a smooth manifold $M$ and $(A, B)$ a pair of  smooth stratified submanifolds of $M$, we say $A$ and 
 $B$ are transverse, denoted $A\pitchfork B$, if the strata of $A$ and $B$ are pairwise transverse.
 \end{defn}
 
 If $A$ and $B$ are submanifolds with $C^1$ conic singularities, one proves easily that mutual
  transversality is a generic property. More precisely, if $A$ and $B$
 are not transverse, generically for a diffeomorphism $g$ of $M$, the manifold $g(A)$ is transverse to $B$.
 
 \begin{prop}If $A$ and $B$ are two mutually transverse submanifolds of $M$ with $C^1$ conic singularities,
 then $A\cup B$ is also such a submanifold. 
 \end{prop}

 \proof (Sketch) Let $S_A$ and $S_B$ denote two strata of $A$ and $B$ respectively. Let $L:= S_A\cap S_B$.
 Consider the normal bundle to $S_A$ in $M$ over $L$. Let $S\nu_A$ denote its sphere bundle. And similarly with $B$.
 Then $L$ is a typical stratum of $A\cup B$. The fiber over $x\in L$ of its conic transverse structure is  isomorphic
 to the cone based on the \emph{join} $(S\nu_A)_x*(S\nu_B)_x$. This join, viewed in the unit sphere of the fiber over $x$
 of the normal bundle to $L$ in $M$, consists of the union of geodesic segments from
 $(S\nu_A)_x$ to $(S\nu_B)_x$.\bull

\section{About conic singularities in Morse theory} \label{conic-sing}
In 1992, I have been asked by Jean-Michel Bismut how is the closure of the stable/unstable manifolds 
of a gradient vector field $X$ of a Morse function $f$ on a compact 
$n$-dimensional manifold $M$. The answer is simpler under the assumption 
the gradient is \emph{simple} near all critical points of $f$ in the following sense.
 
 \begin{defn} The vector field $X$ is said to be a \emph{simple descending gradient}\,\footnote{ For cell attachings 
  it would be uncomfortable to work with an ascending gradient.} of $f$ near $a\in {\rm crit}f$ if 
the   $X=-\nabla f$ in the Euclidean metric
of Morse coordinates $x=(x_1, \ldots, x_n)$ where $f$ reads 
\begin{equation} \label{eq-simple}
 f(x)= f(a)+ \frac 12 (-x_1^2-\cdots-x_k^2 +x_{k+1}^2+\cdots+x_n^2).
\end{equation}
Here $k$ is the Morse index of $f$ in $a$. The set of critical points of index $k$ of $f$ will be denoted 
${\rm  crit}_kf$.
\end{defn} 
In that case, the local stable (resp. unstable) manifold $W^s_{\rm loc}(a,X)$ (resp. $W^u_{\rm loc}(a,X)$)
are radially foliated by the local orbits of $X$.

By Smale \cite{smale},  the gradient $X$ is generically Morse-Smale meaning that the stable and unstable manifolds
of all critical points are pairwise transverse.\footnote{ Smale did not assume $X$ is simple. But the same proof yields 
the genericity of Morse-Smale gradients among the simple gradients.} We have proved the next proposition in \cite{bz}.
It is stated by using Definition \ref{conic} of $C^1$ conic singularities.

\begin{prop}\label{closure} If $f$ is a Morse function on $M$ and $X$ is a Morse-Smale gradient which is simple near the 
critical points of $f$, then the closure of the invariant manifolds $W^{u/s}(a,X)$, $a\in{\rm crit} f$, are submanifolds
 with $C^1$ \emph{conic} singularities.
\end{prop} 

Equivalently, on may state the following: \emph{Under the same assumptions the union
 \begin{equation}
 \Si:=\mathop{\cup}_{{\rm index}(a)<n} W^u(a, X)
 \end{equation} 
 is a submanifold of $M$ with $C^1$ conic singularities.\\
}


\begin{remarques}\label{rem-conic} 
{\rm ${}$

\nd {\bf 1.} In the setting of  Proposition \ref{closure} (2nd form), the stratum $\Si_j$  is made of the union of the unstable manifolds associated to the critical points of Morse index $j$.  In that case,
the regularity of the local trivializations by radial cones is not greater than 1 in general. This regularity issue 
is discussed in \cite[Appendice]{coursX}.

\nd {\bf 2.} As observed in Section 2 of \cite{ab-laud}, Proposition \ref{closure}  
extends to the case where $M$ has a non-empty 
boundary and the data are generic with respect to the boundary.\footnote{ In this setting of non-empty boundary,
the stable/unstable manifolds are explained in \cite{lauden1}.}
}
\end{remarques}
The next corollary---a statement anticipated by Ren\'e Thom \cite{thom}---was given in \cite{bz} just as a remark with a hint for proving; indeed, it was not the main goal
of that appendix. At that time, I was seriously criticized about this remark. Among those who were confident in the 
statement, some told me that difficult steps were missing in my suggested proof. Here, I would like 
to give a formal proof for closing the discussion. 

I should say that this statement was first proved by G. Kalmbach \cite{kalmbach} in the same setting---a reference that I was ignoring up to reading the next reference.
Very recently, the cell decomposition has been obtained 
 by A. Abbondandolo \& P. Majer \cite{abbondan} for every gradient-like vector fields; beside the mathematical part,
this paper contains an interesting history of the topic.

\begin{cor}\label{cell-dec}
In the setting of Proposition \ref{closure}, the partition of $M$ made of the  unstable manifolds $W^u(a,X)$,
 $a\in {\rm crit} f$, is induced by some $CW$-complex structure on $M$ in the sense that
  its ``open" cells are the unstable manifolds 
 of the zeroes of $X$.
\end{cor}

What is missing in the information given by the unstable manifolds is the way where each \emph{closed} cell is 
attached to the lower dimensional skeleton. A priori, it is not immediately clear that the closure of $W^u(a,X)$ is the image
of a closed ball and the attaching map does not seem to be given by the dynamics. 

Before proving Corollary \ref{cell-dec}  we need some refreshing about Proposition \ref{closure} (see the next two subsections).  

\begin{rien}\label{magic} {\sc Magic of the simple Morse coordinates.} {\rm 
Since maxima and minima are not interesting for the present discussion, 
we limit ourselves to the case of the 
$n$-dimensional Morse Model of index $k, \ 0<k<n$, denoted $M(n,k)$. It is obtained by rotating the 2-dimensional 
model of index 1 around the $x$-axis\,\footnote{ In general, this axis is not 1-dimensional! and similarly for the $y$-axis.}
and the $y$-axis in order to create stable and unstable manifolds of dimension $n-k$ and $k$ respectively. 

There is a top 
 boundary denoted $\p^+M(n,k)$ which is a trivial $k$-dimensional disc bundle over some $(n-k-1)$-dimensional Euclidean sphere  $S^+$  centered at $O$ in 
 the local  stable manifold $W^s_{loc}(O,X)$. 
 The fibers are discs parallel to the unstable manifold $W^u_{loc}(O,X)$---for brevity
 the radius of $S^+$ is not specified. 
  Similarly the bottom boundary
 $\p^-M(n,k)$ is a trival $(n-k)$-dimensional disc bundle over a Euclidean $(k-1)$-sphere $S^-$ in $W^u_{loc}(O,X)$;
  its fibers are parallel to $W^s_{loc}(O,X)$.  
 The core of the top (resp. the bottom) boundary is named the \emph{co-sphere}  
 (resp. the \emph{attaching sphere})---see Figure \ref{morse}.

Let $\Si^+\subset \p^+M(n,k)$ be a $(k+s)$-dimensional submanifold with $C^1$ conic singularities and let $F\Si^+$
be the closure in $M(n,k)$ of the set of trajectories starting from $\Si^+$. Suppose $\Si^+$ transverse to the co-sphere 
and set $K:= \Si^+\cap S^+$. Let $cK$ be the radial cone of $K$ in $\{0\}\times D^{n-k}$. 
For every $x\in S^-$ the fiber $D_x$ of $\p^-M(n,k)$ is provided with a cone $cK_x$ by translating over $x$
some homothetic of $cK$.

\begin{center}\label{morse}
 \begin{figure}[h]
 \includegraphics[scale =.6] {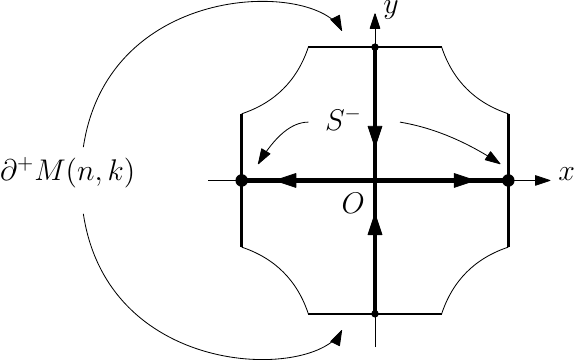}  
 \caption{\small Here, $n=2$, $k=1$.}
\end{figure}
\end{center}

By transversality, the projection of $\Si^+$ to the unstable manifold is of maximal rank near $K$; then,
by flowing, the projection to the factor $D^k\times \{0\}$ of first $k$ coordinates is also onto. So, $F\Si^+$ contains
 $(D^k\times \{0\})\cup( \{0\}\times cK)$ and, for every $x\in S^-$, the disc
$D_x$ is transverse to $F\Si^+$. 

By a beautiful property of the polar coordinates---here, one should speak of multispherical coordinates---the following 
holds.}
\medskip

\nd{\sc Magic Fact.} Whatever the angle along $K$ between 
 $\Si^+$ and the co-sphere,  for every $x\in S^-$, the radial cone $cK_x\cong cK$ is the tangent cone 
to $F\Si^+\pitchfork D_x$.
\medskip

{\rm This is the main flavour of Proposition \ref{closure}. Of course, the same propagates along the unstable manifold
for every choice among 
its tubular neighborhoods which are fine enough.}
\end{rien}

\nd \begin{rien} {\sc Some choices.}\label{choices}
 {\rm By induction on the dimension, one truncates the unstable 
manifolds 
in order to create for every $k$ 
some compact domains $\u\Si_k\subset \Si_k$ 
 and some tubular 
neighborhoods $T_k$ over $\u\Si_k$, sufficiently fine for still having the magic conic property  and fulfilling a few more conditions:
\begin{enumerate}
 \item {\it $T_0$ is a finite union of balls where $X$ is radial.}
 \item {\it  The restriction of $T_{k}$ 
to the boundary $\p\u\Si_k$ is covered by the interior of $T_0^{k-1}:= T_0\cup\cdots\cup T_{k-1}$.}
 \item  {\it Let $W_k$ be a collar neighborhood of $\p\u\Si_k$ in $\u\Si_k$ and let $E_k$ be the restriction of $T_k$ to
$W_k$. Then, if $E_k$ meets $T_j$ and $j < k$, we have $E_k \subset {\rm int\,} T_j$. Moreover, in that case, each fiber of $E_k$  
over $W_k$ is contained in some fiber of this $T_j$.}
\end{enumerate}
The third item adds some requirement to the disc bundle over $\u\Si_k$. But, as it is said 
in the very end of Subsection \ref{magic}, the conic structure is not affected by this new requirement.
}
\end{rien}
  \begin{center}\label{fig-pref}
 \begin{figure}[h]
 \includegraphics[scale =.6] {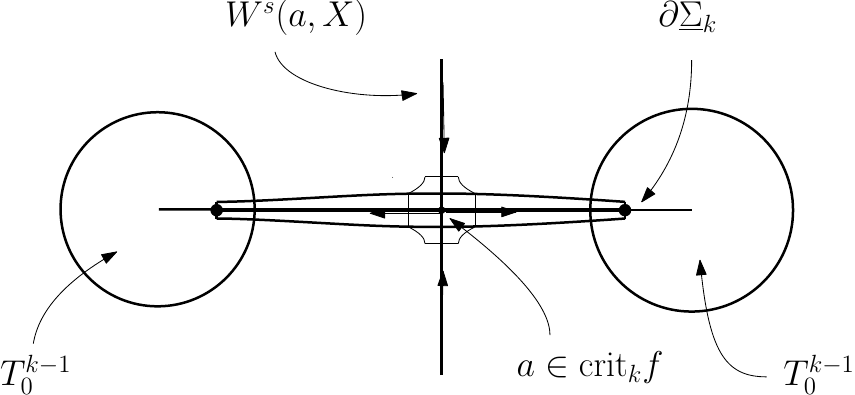}  
 \caption{\small Here, $a$ is a critical point of index $k$. Its unstable manifold  is schemed horizontally. The two balls 
 represent parts of the chosen neighborhood of the $(k-1)$-skeleton.}
 
 \end{figure} 
 \end{center}

\begin{rien} {\sc Proof of Corollary \ref{cell-dec}.} {\rm We are going to prove,  by an induction on $k$
 from $k=0$ to $n$,
that the union of unstable manifolds of dimension not greater than $k$ is a $CW$-complex 
whose open cells are these unstable manifolds.
For $k=0$, there is nothing to do: one starts with a finite collection of points. Assume the claim holds for 
the union of strata of dimension less than $k$.

Let $p_k: T_k\to \u\Si_k$ denote the projection of the bundle structure. The proof may be made cell by cell.
For $a\in {\rm crit}_kf$, we denote $\u W^u(a,X)$ the intersection $W^u(a,X)\cap\u\Si_k$. 
Choose an identification of $i_{op}: W^u(a,X)\to {\rm int} D^k$  
with the open   $k$-ball sending the orbits of $X$ to the Euclidean rays. And choose a similar identification
$i_c: \u W^u(a,X) \to D^k$ with the closed $k$-ball;  these two are supposed 
to coincide near $a$. The composition $i_c^{-1}\circ i_{op}$ yields an \emph{abstract} compactification of $W^u(a,X)$
as a closed $k$-cell.

If we are able to describe a continuous map $f_a$ from $\p\u W^u(a,X)$ to the $(k-1)$-skeleton $\Si_0^{k-1}$
whose image is the frontier of  $W^u(a,X)$,  
then $f_a$ is one (among many others) desired cell attachment. 

If $z\in\p\u W^u(a,X) $ is covered by only one tube $T_j$, $j<k$, then $f_a(z):= p_j(z)$ looks to be an easy recipe for 
cell attaching. But if $z$ is covered by $T_{j}$ and $T_{i}$,  $i<j$, we need to ``average" the two possible 
projections; here, item (3) of Subsection \ref{choices} is used. Concretely, one proceeds as follows.

One chooses some continuous function $\la_j: W_j\to [0,1]$ with value $0$ on $\p\u\Si_j$ and $1$
on the opposite side of the collar; one extends  $\la_j$ to $E_j$ by pre-composing  with $p_j$.
By the conic structure transverse to $\Si_i$ there is an arc $\ga_z:[0,1]\to\Si_i$, depending continuously on $z$
which joins  $x:=p_i(z)$ to $y:=p_j(z)$ in the fiber $T_{i,x}$ and 
$\ga_z(s)$ runs in $\Si_j$ when $s>0$.

Therefore, the following ``averaged" attachment is convenient in the above setting of three strata 
$(\Si_i,\Si_j, \Si_k)$, ranked by increasing dimension and pairwise adhering one to the next:
\begin{equation}\label{3strata}
f_a(z)=\left\{
\begin{array}{ll}
p_j(z) &\text{if } z\in (T_j\smallsetminus E_j)\\
\ga_z(\la_j(z)) & \text{if } z\in E_j\\
p_i(z) & \text{if not.}
\end{array}
\right.
\end{equation}

In general, it is a sequence $(\Si_{i_1}, \Si_{i_2},\ldots, \Si_{i_r}, \Si_k)$ that we have to consider, where 
each stratum adheres to the next one. In that case, there are several projections 
of $z\in \p\u W(a,X)$ to $p_{i_\ell}(z)\in \Si_{i_\ell}, \ 1\leq \ell\leq r$; there are also several 
paths $\ga_{z,\ell}$ joining $p_{i_{\ell-1}}(z)$ to $p_{i_\ell}(z)$ and running in $\Si_{i_\ell}$ at positive time. Then 
the attaching map $f_a(z)$ reads
\begin{equation}\label{r+1}
f_a(z)=\left\{
\begin{array}{ll}
p_{i_r}(z) &\text{if } z\in (T_{i_r}\smallsetminus E_{i_r})\\
\ga_{z,r}(\la_{i_r}(z)) & \text{if } z\in E_{i_r}\\
\cdots & \cdots\\
\ga_{z,\ell}(\la_{i_\ell}(z)) & \text{if } z\in (E_{i_\ell}\smallsetminus E_{i_{\ell-1}})\\
\cdots & \text{and so on} \\
p_{i_1}(z) & \text{if not.}
\end{array}
\right.
\end{equation}
So, the attaching map of $\u W(a,X)$ is now completely described and Corollary \ref{cell-dec} is proved. \bull
}
\end{rien}

\section{Immediate transversality}\label{imm-trans}

Even in their advanced version, for instance in jet spaces \cite{jet-space} or in a more abstract context \cite{thom-trans},
Thom's transversality theorems just state approximation results, possibly with fine topologies 
on the considered functional spaces. 

These transversality theorems were not sufficient for getting $A_\infty$-structures on
 Morse complexes in the line of K. Fukaya's program \cite{fukaya}. 
  The enhanced version of transversality, we named \emph{immediate transversality}, appeared to be the right 
  tool---possibly among others---to complete this program as we did  in \cite{ab-laud}. 
  This seems to be a new topic.
 
 In the present paper,  we slightly generalize the setup compared with \cite{ab-laud}, but we are far from knowing
 what is the right setting for immediate transversality.
 
 \begin{defns} ${}$\label{immediate}
 
 $(1)$ Let $f: A\to M$  be a smooth map and $S$ be a smooth submanifold in $M$. The ambient 
 isotopy $(\vp^t)_{t\in [0,1]}$ of $M$, from $Id_M$, is said to make $f$ \emph{immediately transverse} to $S$
 if there exists $\ep>0$ such that
 \begin{equation}
 (\vp^t\circ f) \pitchfork S \text{ {\rm for every}  } t\in (0,\ep).
  \end{equation} 
  
  $(2)$ When $\vp^t$ in (1) is the flow of an autonomous vector field one speaks of 
  \emph{immediate transversality by flow}.
 \end{defns}
 
 It is worth noting that, in case of Definition \ref{immediate}(2), one gets a remarkable and very rare property, namely 
 the transversality to every element of 
  a given one-parameter family. More precisely, for every $0<t_1<t_2<\ep$ the following holds:
 \begin{equation}\label{trans-path}
(\vp^{t_2} \circ f)\pitchfork \vp^{t} (S) \text{ for every } t\in [0,t_1].
 \end{equation}
 For proving the above formula, one just applies the one-parameter group formula of flow and the fact that an ambient diffeomorphism 
 preserves transversality. 
 
 \begin{rien} {\sc Toy example.} {\rm Given a closed $n$-dimensional smooth manifold $M$ endowed 
 with a closed smooth submanifold $S$, one can consider a closed tubular neighborhood $N$ of $S$ in $M$, 
 with a projection $p:N\to S$,
 and a smooth section $u: S\to N$ vanishing transversely. Such a section exists by the very first transversality theorem of Thom \cite[Th\'eor\`eme I.5]{elementary}. 
 
 Regarding $(N,p,S)$ as a linear disc bundle one is allowed, for $x\in S$, to think of $u(x)$ as a vector
 in the vector space generated by the disc $N_x$, the fiber over $x$. Moreover, thanks to the underlying affine structure, 
 one can translate $u(x)$ by parallelism at every point of $N_x$. So, we have a vector field, denoted $\vec u$, 
 over $N$ that extends arbitrarily
 to the complement of $N$ in $M$.
 
 The germ of its flow is a translation flow: for $x\in S$ and every small enough
  positive $t$, we have $\vec u^{\,t}(x)=x+t u(x)$. Since
 $u$ is a section of $T$ transverse to the zero-section, the flow generated by $\vec u$ makes
  $S$ immediately transverse
 to its initial position. We shall see that the case of submanifolds with conic singularities is less elementary.
 }
 \end{rien} 
 
 \begin{thm} \label{imm-trans-conic}
 Let $\Si$ be a compact submanifold with $C^1$ conic singularities in an $n$-dimensional smooth manifold $M$.
  Then there exists 
 a ``large'' family of flows $(\vp^t)$ on $M$ making $\Si$ immediatly transverse to its initial 
 position.\footnote{ Here, a family $\mathcal F$ of maps $C\to E$,
parametrized by a manifold $P$, is said to be \emph{large} if the corresponding map 
$F: P\times C\to E$ is a submersion at every point of its domain. One may extend this definition to families of flows
by demanding this property at every positive time.}

 \end{thm}

The main example we have in mind is given by the union of 
unstable manifolds of a simple descending gradient
of a Morse function excluding the $n$-dimensional unstable manifolds.
The difference with \cite{ab-laud} is that no particular assumptions are made on the  topology of  strata in $\Si$.\footnote{ In \cite{ab-laud} the strata are disjoint unions of Euclidean spaces.}\\

\proof It will be done in three steps. In Step I, one looks at the neighbourhood of one vertex of $\Si$, 
that is, the cone based on 
a stratified submanifold in the $(n-1)$-sphere. In Step II, one looks at a $C^1$ bundle over a compact base 
whose fibers are stratified cones. 
In Step III, one glues together the different pieces making  the given submanifold with
$C^1$ conic singularities.\\

\nd{\sc Step I: a unique radial cone.} {\rm One may suppose the ambient manifold is $\R^n$
provided with a cone $C$ whose vertex is the origin $O$ and which is 
based on a compact stratified submanifold $L$ (as \emph{link})
of the $(n-1)$-dimensional sphere $S^{n-1}$ with $C^1$ conic singularities.
The strata of $C$, apart from the origin, are the cones on strata of $L$, punctured at $O$. 

\begin{prop}\label{step1} 
Generically for $\vec u\in \R^n$, the translation flow generated by $\vec u$ makes $C$ immediately transverse
to $C$.
\end{prop}

\proof One first looks at the space of secants
$Sec(S_1, S_2)$ for each pair $(S_1, S_2)$ of strata from $C$. We choose to parametrize $Sec(S_1, S_2)$
as a submanifold of $S_1\times \vec\R^n$:
\begin{equation}Sec(S_1, S_2)=\{(x,\vec u)\in S_1\times \vec\R^n \mid x+\vec u \in S_2\}
\end{equation} 
Without loss of generality we may suppose that $S_i$, $i=1,2$, has a regular equation $f_i=0$. Therefore,
one checks that $Sec(S_1, S_2)$ is defined by the system
\begin{equation}
\left\{
\begin{array}{l}
f_1(x)=0\\
f_2(x+\vec u)=0
\end{array}
\right.
\end{equation}
whose linearized system reads
\begin{equation}
\left\{\begin{array}{l}
df_1(x)\,\de x=0\\
df_2(x+\vec u)(\de x+\de\vec u)=0.
\end{array}
\right.
\end{equation}
}
As $x$ and $\vec u$ are independent variables, it follows that this system is of maximal rank. So, $Sec(S_1, S_2)$
is smooth. Then, Sard's theorem applies: the projection 
\begin{equation}
\begin{array}{c}
\pi: Sec(S_1, S_2)\to \vec\R^n\\
(x,\vec u)\mapsto \vec u
\end{array}
\end{equation}
has a dense set---in the sense of Baire which is now called \emph{generic}---of regular values.\\

\nd {\sc Claim: Coplanarity criterium}. {\it The vector $\vec u\in\vec\R^n$ is a critical value of $\pi$ if and only if 
the tangent spaces $T_xS_1$ and $T_{x+\vec u} S_2$ are \emph{coplanar} that is, are contained 
in a common hyperplane.\footnote{ If $H$ is a hyperplane tangent to $S_1$ at $x$ it must contain the corresponding
 generator of the cone $S_1$, and hence the origin.}
}\\

\nd {\sc Proof $\Leftarrow$.} Suppose $L$ is a linear form on $\vec\R^n$ which vanishes   on 
the two kernels $\ker df_1(x)$ and $\ker df_2(x+\vec u)$. Then $L(\de x)=0$ and $L(\de x+\de\vec u)=0$
for every $(\de x,\de\vec u)$ tangent to $Sec(S_1, S_2)$. Therefore, $L(\de\vec u)=0$ and the image of $d\pi_(x,\vec u)
$ is contained in the kernel of $L$. So, $\pi$ is not a submersion near $(x,\vec u)\in Sec(S_1, S_2)$.\\

\nd $\Rightarrow$. By contraposition, we have to prove that if the two considered kernels are not coplanar 
then, for every $\de \vec u\in \vec\R^n$, the system 
\begin{equation}
\left\{
\begin{array}{ll}
df_1(x)\,\de x&=0\\
df_2(x+\vec u)\,\de x&= -df_2(x+\vec u)\,\de\vec u
\end{array}\right.
\end{equation}
has a solution with  $\de x$ as unknown. This is an exercise of Algebra which we leave to the reader. \bull

This claim still holds when considering all pairs of strata; we just have to define the tangent space $T_xC$
as the tangent space to the unique stratum which contains $x$. Indeed, it sufficies to intersect 
finitely many residual sets since there are finitely many strata numbered by their dimension.

Having this claim in hand, we conclude that, generically, $\vec u$ translates the cone $C$ to a cone 
which is transverse to $C$, that is 
$(C+\vec u)\pitchfork C$. But elementary geometry shows that the critical values of $\pi$ is itself a cone.
As a conclusion, the translation flow generated by a generic vector $\vec u$ is of immediate transversality
with respect to $C$ which what is desired for Step I.\bull

We have not used the compactness of the basis of $C$. This additional assumption allows us to state 
that the residual set $R$ of generic vectors
in question is the complementary of a cone with a compact basis.  
So, $R$
is open in $\vec \R^n$. This is proved in a more general setting in Proposition \ref{openness} (see also \cite[Appendix B]{ab-laud}.) \\

\nd {\sc Step II: A $C^1$ cone bundle over a compact base.} {\rm
Let $(E,C)$ be the pair made of an $(n-k)$-disc bundle 
and a $C^1$ conic subbundle;
 their  common base space $B$ is a  compact $k$-dimensional smooth manifold. 
 Up to taking a smaller tube, by Proposition \ref{linearization} we may assume $C$ is a linear conic subbundle
 without loss of generality.

Let $\mathcal U=\left\{U_\al\right\}_{\al\in I}$ be a finite 
open covering of $B$
with the following properties where the index set is  
$I= \{1, \ldots, \ell\}$:
\begin{enumerate}
\item For every $\al\in I$ the pair  $(E,C)$ has a \emph{linear}  trivialization $\phi_\al$ over the open 
set $U_\al$.   
\item Each $U_\al$ is endowed with a compactly supported smooth function $\rho_\al:U_\al\to [0,1]$ such that,  
naming $V_\al$
the interior of the level set $\rho_\al^{-1}(1)$, the family $\{V_\al\}_{\al\in I}$ still covers $B$.\\
\end{enumerate}

\nd {\sc Step II-1: $C$ is a linear conic subbundle of $E$.}

Then one considers the vector space  $G$ of $\ell$-tuples 
 ${\bf v}:= (\vec v_1, ... , \vec v_\ell)$ of vectors in $\R^{n-k}$, that is, 
$ G:=\prod_{i=1}^\ell\vec\R^{n-k}$. 
Such a ${\bf v}$ acts on $E$ by translation in each fiber $E_x$, $x\in B$.
More precisely, the vector  ${\bf v}=(\vec v_\al)_{\al\in I}$ acts on $z\in E_x$ by the formula
 \begin{equation}
 {\bf v\,\cdot\,}z=(\vec v_1, \vec v_2, \cdots, \vec v_\ell)\cdot z = z+\sum_\al \phi_{\al,x}
 \bigl(\rho_\al(x)\vec v_\al\bigr)= : z+\si_{\bf v}(x).
 \end{equation}
This formula is well defined since $\phi_\al$ is linear---up to replacing the disc bundle $E$   by the
spanned vector bundle $span(E)$---and two translations commute. An element ${\bf v}\in G$ 
may be thought of as a translation field in $E$ over the base $B$ that we denote $\si_{\bf v}$.

For every $x\in V_\al$ and $z\in E_x$, this action is submersive at $z$. Indeed,  this action is already 
submersive if the entries of ${\bf v}$  are null except the entry of index $\al$. So, we have 
a map $\psi: G\times C\to span(E)$, piecewise smooth and submersive onto $E$
at every point of the source. The inverse image $\psi^{-1}(C)$
is  stratified and Sard's theorem applies to the projection $\pi: \psi^{-1}(C)\to G$. So, the set of critical values of 
$\pi$ is meager. By characterizing these critical values
similarly to Step I, one concludes
there is a residual $R\subset G$ such that for every ${\bf v} \in R$ we have 
\begin{equation}\label{trans-cone-fam1}
({\bf v\,\cdot}\,C )\pitchfork C,\text{ that is } (C+\si_{\bf v}) \pitchfork C,
\end{equation} 
where $\si_{\bf v}$ stands for the section $x\mapsto {\bf v}(x)$ and the sum is applied  fiberwise.

 Now, if $\si$ is an arbitrary other smooth section of $span(E)$ the family 
 $\{{\bf v\,\cdot\,}\si\}_{{\bf v}\in G}$ is a submersive family onto $E$.
   So, by the same reasoning as in the case $\si=0$ we get
   a dense set of smooth sections  $\tilde\si$ of  $span(E)$  such that 
     \begin{equation}\label{sec-trans}
     (C+\tilde\si)\pitchfork C. \footnote{ This way of reasoning by using 
  finite dimensional submersive families 
  is essentially what is done by Thom in  \cite{thom-trans}.}
     \end{equation} 
  Moreover, due to the radial conic feature of $C$, such a section generates a flow of immediate transversality to $C$;
  more precisely, 
 \begin{equation} \label{t-sec}
   (C+t\tilde\si)\pitchfork C\text{ for every } t>0.
  \end{equation}
 This may be checked by showing that the opposite of (\ref{t-sec}) is a homogeneous condition on $\si$.\footnote{
 At that aim it is convenient to look at the ``completed cone bundle'', that is, in $span(E)$ instead of the disc bundle.}
    So, we have:

  \begin{prop} \label{family-linear}
  Let $C$ be a linear cone bundle in a disc bundle $E$ over $B$.
    Then, there is a dense set of smooth sections
    $\si $ of $span(E)$ which generates a translation flow on $span(E)$ of immediate transversality to $C$.\\
 \end{prop}
 
 \begin{remarque}\label{rem-trans}
 {\rm  
 The transversality stated in (\ref{trans-cone-fam1}) does not mean that it holds fiberwise. 
 The latter only holds in generic fibers. In their complement, that is over points of the bifurcation locus,
  transversality is provided by the image of the linear derivative operator
     $D\si_{\bf v}$. More precisely, the section $\si_{\bf v}$, viewed as a translation field over $B$, defines a
     \emph{gauge transformation} $\Theta: span(E)\to span(E)$ over $Id_B$. At $z\in E_x$,
     its derivative $D\Theta: T_zE\to T_{z+\bf v(x)}span(E)$ is a linear map over $Id\vert_{T_xB}$.
     If $x$ belongs to the bifurcation locus, property (\ref{trans-cone-fam1}) states as follows:  For every 
  $z\in C_x$ and every hyperplane $H\subset T_{z+\bf v(x)}E$ which is tangent to the fiber $C_x$ at $z$
  and to  $C$ at $z+{\bf v}(x)$, 
  then
  \begin{equation}\label{trans-cone-fam2}
  D\Theta ( T_zC )\text{ is  transverse to }H.
  \end{equation}
  }
    \end{remarque}
 
 \bigskip

 \nd {\sc Step II-2: Openness property.}

\begin{prop}\label{openness}
 If $C\to B$ is a linear conic subbundle in $E\to B$ over a compact base $B$, the set of smooth sections
 $\si$ of $span(E)$ such that $C+\si$ is transverse to $C$ is  open  in the set of smooth 
 sections of $span(E)$ for the smooth topology. 
\end{prop}

By (\ref{t-sec}) we already know that $(C+\si)\pitchfork C$ implies that $\si$ generates a translation flow of immediate
transversality to $C$.\\

\proof   Let $\si_0$ be a smooth section of $span(E)$ having the demanded  transversality property.
Since the base space of the fibration is compact the question is local over $B$. There are two types of points $x\in B$.

1) {\it Generic points}: that is, $\si_0(x)$ maps $C_x$ transversely to $C_x$ in $E_x$. In that case, the same occurs
  in the nearby fibers. Therefore, for every hyperplane $H\subset T_{z+\si_0(x)}E_x$ tangent to $C_x$ at 
 $z+\si_0(x)$, the angle between $C_x+\si_0(x)$ and $H$---considering each stratum 
  separately---is locally positively bounded from below; here we use the compactness of the considered 
  Grassmanian.
  
  Since this angle does not approach to zero near $z$,
  any small enough perturbation of $\si_0$ in the $C^1$ topology has the same property. That would prove the 
  desired openness over a neighborhood of such a point $z$ if the stratum of $C$ containing $z$ were compact, 
  which is not.
  
  For circumventing this difficulty we are going to use the fact that  
  the spherical basis $S_x\subset \mathbf S^{n-k-1}$ of the cone $C_x$ has itself $C^1$ conic singularities;
  that compensates for non-compactness of the strata of $C_x$.

  Indeed, let   $S_0\subset C_x$
  be a stratum transverse to $H$.
  Then the angle between $S_0$ and $H$ controls, at  every point close enough to $S_0$, the angle of $S$ with $H$ 
   for every stratum $S\subset C_x$ adhering to $S_0$. As a consequence, the reasoning on angles 
  works as if the strata of $C_x$ were compact.\\
  
2) {\it Bifurcation locus.} Now, let us consider a point $x_0\in B$ where  $(C_{x_0}+\si_0(x_0))\pitchfork C_{x_0}$ is not fulfilled. Nevertheless, by assumption on 
$\si_0$ for every $z\in C_{x_0}$ we have 
\begin{equation}
 D\Theta(T_zC)\pitchfork H
\end{equation}
for every hyperplan   $H\subset T_{z+\si_0(x_0)}E$ tangent to $C$ at $z+\si_0(x_0)$ and to $C_x$ at $z$; here, $\Theta$
stands for the affine automorphism of $span(E)$ defines by the translations $\si_0(x)$, $x$ close to $x_0$.
(compare with Remark \ref{rem-trans}.)   
One recalls that $T_zC$ means the tangent
space to the stratum which contains $z$. Again we get an angle  which has a positive infimum when
$z$ runs over $C_{x_0}$. This angle varies continuously under $C^1$ perturbation,  
whatever the nature of $x_0$, generic or not, in the perturbed translation.\bull

As a consequence of this openness property, we can state and prove a relative version of Proposition \ref{family-linear}
whose statement is now enhanced in the following form.

\begin{cor} \label{relative}
 In the above setting, let $\si_0$ be a germ of generic section of $span(E)$ over a collar neighborhood 
of the boundary $\p B$. Then, as a germ, $\si_0$ extends to a global section $\si$ of $span(E)$ 
which generates a translation flow of immediate transversality to $C$. Moreover, the set of such extensions
is open and dense in the space of all extensions of the given germ.
\end{cor}

\proof For the proof, it is easier to realize the germ $\si_0$ as a section over some collar neighborhood $W$ of $\p B$.
Let $\tilde \si_0$ be an arbitrary extension of this realization of the germ, still denoted by
 $\si_0$, to a global section of $span(E)$. Let $\la: W\to [0,1]$ be a  smooth function vanishing on a neighborhood 
 of $\p B$ and equal 
 to $1$ near the interior boundary component of $W$.
   
By Proposition \ref{family-linear}, there exists a global section $\tilde\si$ which generates a translation flow
of immediate transversality to $C$ and may be chosen arbitrarily close to $\si$. 
By Proposition \ref{openness}, if this approximation 
is small enough in the $C^1$ topology, the section 
$$\si_{W}:=\si_0+ \la(\tilde\si_{\vert W}-\si_0)$$
 is still a generic section
over $W$. Since genericness is a local property, the gluing of 
$\si_{W}$ and $\tilde\si_{\vert\overline{B\smallsetminus W}}$
yields the desired section (by the choice of $\la$, the germ of $\si_W$ along $\p B$ is the germ $\si_0$.)\break \bull

\nd {\sc Step III: Gluing of tubes.} Recall the compact submanifold $\Si$ in $M$ with $C^1$ conic singularities.
By Proposition \ref{linearization} we may assume the conic structure induced by $\Si$ in the tube chosen around 
each stratum is linear in the sense of Definition \ref{lin}.
Let $S_k$ and $S_j$ be two connected components respectively of $j$- and  
$k$-dimensional strata in $\Si$. 
By the choices made in subsection \ref{choices}, we have compact domains $\u S_j$ and $\u S_k$; and also an
$(n-j)$-disc bundle $N_j$ over $\u S_j$ (resp. an $(n-k)$-disc bundle over $\u S_k$.)

In what follows, $S_j$ lies in the closure of $S_k$. By conditions (2) and (3) in the mentioned choices,
the collar neighborhood $W_k$ of 
$\p\u S_k$ in $\u S_k$ meets $N_j$ and avoids its sphere bundle $SN_j$; and similarly for the disc-bundle $E_k$
induced by $N_k$ over $W_k$.
Finally, each fiber of $E_k$
generates an affine $(n-k)$-subspace in a fiber of $N_j$.\\

\nd{\sc Step III-1: Reduced translation flow and reducing process.} 
Let $u$ be a translation field in $N_j$. For every fiber $N_{j,z}$,
every point $x\in W_k\cap N_{j,z}$ and every $y\in N_{k,x}$ we have a splitting of tangent spaces
\begin{equation}
T_y (N_{j,z}) = T_x( W_k\cap N_{j,z})\oplus T_y (N_{k,x})
\end{equation}
Here, the tangent space $T_x( W_k\cap N_{j,z})$ is carried  up to point $y$ by parallelism with respect to the affine 
structure of $N_{j,z}$. 
This splitting decomposes the translation vector $u(z)=u(y)$ into horizontal and vertical component
at $x$, namely:
\begin{equation}
u(z)= u_\frak h^k(x)\oplus u^k_\frak v(y)
\end{equation}
This splitting is independent of $y$ along the fiber $N_{k,x}$. The vertical component $x\mapsto u^k_\frak v(x)$
is termed the \emph{reduction of }$u$ to $N_k$.

The bundle structure of $\Si\cap N_j$ over $\u S_j$ allows us to choose the collar $W_k$ so that 
its intersection with $N_{j,z}$ is a collar of $\p\u S_k  \cap N_{j,z}$ in $\u S_k\cap N_{j,z}$. 
One chooses a smooth function
$\mu: W_k\to [0,1]$ equal to 1 near $\p\u S_k$ and equal to 0 near the opposite boundary component of the collar. Then
$\mu$ is lifted to $E_k$ as a constant function in each fiber $N_{k,x}$ over $W_k$.
The \emph{balanced reducing process} consists of replacing the constant vector field $u(z)$ on $E_{k,z}$
by the vector field 
\begin{equation}
u^k_\mu(x):= \mu(x)u^k_\frak h(x)+u^k_\frak v(x)
\end{equation}
which is termed the \emph{balanced reduction of $u$ to $T_k$}. Note $u^k_\mu$ is equal to $u$
over $\p\u S_k$ and to the reduction of $u$ to $T_k$ over the opposite boundary component of the collar $W_k$.

\begin{remarque} \label{seq-adherence}
{\rm In general, at some points of $\p\u S_k$ one has to consider not only the inclusion 
$S_j\hookrightarrow cl(S_k)$ into the closure of $S_k$ but a sequence 
$S_j= S_{j_1}, S_{j_2}, ...S_{j_p}=S_k $ where each term (that is a connected component of stratum) lies in the 
closure of the next term---compare with formula (\ref{r+1}). This would require to compose a sequence 
of reducing processes with a lot of notation. As there is no new idea in such details, they will be skipped.}\\
\end{remarque}

\nd{\sc Step III-2} For completing the gluing, we just need the next proposition. The setting is the same as in 
the previous step. 

\begin{prop}\label{gluing} The translation field $u$ on the tube $N_j$ is assumed to generate a flow
 of immediate transversality to $\Si\cap N_j$.
Then we have:
\begin{enumerate}\item The reduction of $u$ to $N_k$ generates a flow of immediate transversality 
to $\Si\cap N_k$.
\item The balanced reduction of $u$ to $N_k$ generates a flow of immediate transversality to $\Si\cap E_k$.
\end{enumerate}
\end{prop}

\nd {\bf Proof.} As it rarely occurs, 
this  issue deals with bi-1-jets\footnote{ This matter is known to be delicate to handle.} of $N_j$, that is pairs of tangent planes to $\Si\cap N_j$ and is local in nature. So, it is enough to look at the order-one Taylor expansion of the considered vector fields.

Over the neighborhood of the considered $x\in W_k$, over there $E_k$
fibers not only over $W_k$ but also over some domain of $\u S_j$.  So,  
we are allowed to  choose a local trivialization $\varphi$ of $(E_k, \Si\cap E_k)$ 
which extends to a local trivialization $\tilde\varphi$ of $(N_j, \Si\cap N_j)$.

 After such a choice,  the linear disc bundle 
$N_k$ is locally endowed  with two \emph{linear} connections, $h_0$ and $h_1$, seen as plane distributions 
complementary to the fibers.   The distribution $h_0$ splits as a direct sum
 $h_0= h'_0\oplus h''_0$, the first factor $h'_0$ is given by the restriction to $N_k$ of the connection 
$h$ of $N_j$ which is derived
from the local trivialization $\tilde\varphi$
  and the second factor $h''_0$ is parallel to $T_x(\u S_k\cap N_{j,z})$ with respect to the affine structure of $N_{j,z}$
 at every point of the 
fiber $N_{k,x}$. And $h_1= h'_1\oplus h''_1$ where $h'_1=h'_0$ and $h''_1$ is derived from $\varphi$. The difference between these two connections, seen as a vertical deviation,
is measured by a differential 1-form $\om$ on $\u S_k$ valued in the vector space of linear endomorphisms of
$span(N_k)$. 

The translation field $u$ on $N_j$ is assumed to generate a flow of immediate transversality to $\Si\cap N_j$. 
In particular, for every small enough $\ep>0$, there is $\kappa>0$ satisfying the following:
\begin{itemize}
\item[$(*)$] For every $t\in [0,\ep)$,
$z\in \u S_j$,  $x\in W_k\cap N_{j,z}$,  $y \in N_{k,x}$ and for every hyperplane 
$H$ in the tangent space $T_{y+tu(z)}\Si$, then or else $y+tu(z)$ does not belong to $\Si$ or 
the angle $\ga_H(t)$ of $H$ 
with the image of $T_y\Si$ through $T_y(u^t)$ is larger than $\kappa t$, where $T_y(u^t)$ denotes 
the tangent map at $y$ of the translation by $tu$.\footnote{ We recall $T_y\Si$  stands for the tangent space
at $y$ to the stratum of $\Si$ which contains $y$}
\end{itemize}
  So, the constant  $\kappa$ is a velocity in nature. Let $D_zu$ denote the covariant derivative at $z\in \u S_j$ 
  of the section $u$ of $N_j$ with respect to the (local) above-mentioned  linear connection $h$. 
  By linearity of $h$, although $y\in N_{j,z}$ may not belong to the zero section,
it makes sense to derive $u$ at $y$.
  An infinitesimal statement which implies $(*)$ reads as follows: 
  \begin{itemize}
  \item[$(**)$] For every $y\in \Si\cap N_{j,z}$, the subspaces $T_y\Si$ and $D_yu (T_y\Si)$
  are not coplanar in the tangent space $T_yN_j$.\\
  \end{itemize} 
  
\nd{\sc Claim.} \emph{Let $V$ be an open set in $E_k$ and $y$ be a point in $V\cap N_{k,x}$.
If $h_0=h_1$ along $V$, then the statement holds in $V$.}\\

Indeed,   $h_0$ is tangent to $\Si$ as is $h_1 $. Then transversality translates to the vertical component
of $u$. Let $D_x$ be the covariant derivative at $x$ of a section of $N_k$
with respect to the connection $h_1$. 
By an order-one Taylor expansion at $y$, an ``infinitesimal contact'' at $y$, namely,
some hyperplane  $H$ makes $D_y u^k_\frak v(T_y\Si)$ coplanar to $T_y\Si$,
 implies  that some angle $\ga_{H(t)}$ is $o(t)$ as $t$ goes to $0$ and $H(t)$ goes to $H$. This contradicts $(*)$
and proves (1) under the assumption of the claim. If (2) fails it should fail infinitesimally which is impossible by (1).
${}$\bull 

Let $y\in N_{k,x}$ and let $y+tu^k_\frak v$ be the vertically displaced point with a small $t$; suppose both points belong
to $\Si$. The planes $h_1(y)$ and $h_1(y+tu^k_\frak v)$ are both tangent to  $\Si$ but could no longer be parallel
with respect to the $h_0$-parallelism.
Nevertheless, thanks to the 1-form $\om$ which measures the ``difference $h_1-h_0$'', one computes that the angle
between $h_1(y)$ and $h_1(y+tu^k_\frak v)$, the latter being translated to $y$ by a  vertical translation which 
is horizontal with respect to $h_0$\,---whose factor $h''_0$ 
induces the parallelism of the fiber $N_{j,z}$\,---is a $O(t)$.\footnote{ For a given $x\in W_k$, the linear form 
$\om(x)$ evaluated on a vector $\de x\in T_xS_k$ is an endomorphism of $span(N_x)$ which has to be evaluated on  the vector $tu^k_\frak v$.}
Therefore, if $t>0$ is sufficiently small, this angle is negligeable with respect to $\kappa$. So, the reasoning for the claim 
still holds and Proposition \ref{gluing} is proved. \bull

 \medskip
 
\nd{\bf End of the proof of Theorem \ref{imm-trans-conic}.}  Observe that having Proposition \ref{gluing} in
 hand allows one to apply Corollary \ref{relative} for getting 
an extension $\tilde u$ of $u$ to $N_j\cup N_k$ which generates a flow of immediate transversality 
to $\Si\cap(N_j\cup N_k)$. That is what we call a successful gluing! Up to some skipped details 
(Remark \ref{seq-adherence}), the gluing process completes the proof of Theorem 
\ref{imm-trans-conic}: starting from Step I, that is from the  tube $N_0$ with a generic translation $u$, 
one inductively glues the tubes and applies Proposition \ref{gluing}
up to the stratum of maximal dimension which itself offers us a case of toy example. \bull

\end{document}